\documentclass[11pt,letterpaper]{amsart}
\usepackage{amsmath}
\usepackage{amssymb}

\usepackage{epsfig}
\usepackage[all, knot]{xy}
\xyoption{arc}

\newtheorem{thm}{Theorem}[section]
\newtheorem{df}[thm]{Definition}
\newtheorem{prop}[thm]{Proposition}
\newtheorem{cor}[thm]{Corollary}

\newtheorem{lem}[thm]{Lemma}
\newtheorem{ex}[thm]{Example}
\newtheorem{rem}[thm]{Remark}

\newcommand{\Pic}{\operatorname{Pic}}

\newcommand{\Pre}{\operatorname{Preper}}

\newcommand{\pp}{\mathbb{P}}
\newcommand{\af}{\mathbb{A}}

\newcommand{\Rat}{\operatorname{Rat}}

\begin{document}

\title[Height bound for jointly regular families]
{Height bound and preperiodic points
\\for jointly regular families of rational maps}
\author{Chong Gyu Lee}


\date{\today}


\address{Department of Mathematics, Brown University, Providence RI 02912, US}

\email{phiel@math.brown.edu}

\maketitle


\begin{abstract}

    Silverman \cite{S4} proved a height inequality for jointly regular family of rational maps and the author \cite{Le2}
    improved it for jointly regular pairs. In this paper, we provide the same improvement for jointly regular family;
    let $h:\pp^n(\mathbb{Q}) \rightarrow \mathbb{R}$ be the logarithmic absolute height on the projective space, let $r(f)$ be the $D$-ratio of a rational function $f$ which is defined in \cite{Le2} and let $\{f_1, \cdots, f_k\}$ be a finite set of rational maps which is defined over a number field $K$.
    If the intersection of all indeterminacy loci of $f_l$ is empty, then
    \[
    \sum_{l=1}^k \dfrac{1}{\deg f_l} h\bigl(f_l(P) \bigr) > \left( 1+ \dfrac{1}{r} \right) f(P) - C
    \]
    where $r = \displaystyle \max_{l} r(f_l)$.
\end{abstract}


\section{Introduction}

    Let $K$ be a number field and $h:\pp^n(\overline{K}) \rightarrow \mathbb{R}$ be the logarithmic absolute height on the projective space. If $f:\pp^n(\overline{K}) \rightarrow \pp^n(\overline{K})$ is a morphism defined on $K$, then we can make a good estimate of the height $h(P)$ with $h \bigl( f(P) \bigr)$. We can define the degree of given morphism algebraically;
    \begin{df}
    Let $g:V(\overline{K}) \rightarrow W(\overline{K})$ be a rational map. Then, we define the degree of $f$ to be
    \[
    \deg g : = [\mathcal{C}\bigl( V(\overline{K}) \bigr) : g^*\mathcal{C}\bigl( W(\overline{K}) \bigr)]
    \]
    where $\mathcal{C}\bigl( V(\overline{K}) \bigr), \mathcal{C}\bigl( W(\overline{K}) \bigr)$ is the function field on $V(\overline{K})$ and $W(\overline{K})$ respectively.
    \end{df}
    If $f:\pp^n(\overline{K}) \rightarrow \pp^n(\overline{K})$ is a morphism on a projective space, we can find the degree from geometric information;  
    \[
    f^*H = \deg f \cdot H ~\text{in}~\Pic(\pp^n).
    \]
    Then, the functorial property of the Weil height machine will prove the Northcott's theorem. The author refer \cite[Theorem B.3.2]{SH} to the reader for the details of the Weil height machine.
    \begin{thm}[Northcott \cite{N}]
        If $f:\pp^n(\overline{\mathbb{Q}}) \rightarrow \pp^n(\overline{\mathbb{Q}})$ is a morphism defined over a number field $K$, then
        there are two constants $C_1$ and $C_2$, which are independent of point $P$, such that
        \[
            \dfrac{1}{\deg f} h\bigl( f(P) \bigr) +C_1 > h(P)
            > \dfrac{1}{\deg f} h\bigl( f(P) \bigr) - C_2
        \]
        for all $P\in \pp^n(\overline{\mathbb{Q}}) $.
    \end{thm}
    
    If $f$ is not a morphism but a rational map, then the functoriality breaks down; two height functions
     $h_{f^*H}(P)$ and $h_H\bigl( f(P)\bigr)$ are not equivalent. Hence, Northcott's Theorem is not valid for rational maps. (However, we still have $h(P) > \frac{1}{\deg f} h\bigl( f(P) \bigr)+C_2 $ by the triangular inequality. See \cite[Proposition B.7.1]{SH}.)

    Silverman \cite{S4} suggested a way of studying height for rational maps using jointly regular family.
    \begin{df}
    Let $S$ be a finite set of rational maps defined over a number field $K$:
    \[
    S=\{ f_1, \cdots, f_k~|~f_l : \pp^n(\overline{K}) \dashrightarrow \pp^n(\overline{K}) \}
    \]
    and $Z(f)$ be the indeterminacy locus of $f$. We say $S$ is {\em jointly regular} when
    \[\displaystyle \bigcap_{l=1}^k Z(f_l) = \emptyset.\]
    We also say that a finite set of affine morphisms $S'=\{ g_1, \cdots g_k ~|~ g_l : \af^n(K) \rightarrow \af^n(K)\}$ is {\em jointly regular} if corresponding set of rational maps $S= \{f_l~|~ f_l~\text{is the meromorphic extension of}~g_l \in S' \}$ is jointly regular.
    \end{df}
    Then, a jointly regular set will bring an upper bound of $h(P)$;
    \begin{thm}[Silverman, 2006]\label{jointly}
        Let $\{f_1, \cdots, f_k ~|~ f_l : \af^n(K) \rightarrow \af^n(K) \}$ be a jointly regular family of rational maps defined over $K$. Then, there is a constant $C$ satisfying
        \[
        \sum_{l=1}^k \dfrac{1}{\deg f_l} h\bigl( f_l(P) \bigr) > h(P) - C
        \]
        for all $P \in \af^n(\overline{K})$.
    \end{thm}

    In this paper, we will improve Theorem~\ref{jointly};
    \begin{thm}\label{main}
        Let $H$ be a hyperplane of $\pp^n(\overline{K})$, let $S=\{ f_1, \cdots, f_k ~|~ f_l :\af^n(K) \rightarrow \af^n(K) \}$ be a jointly regular family of affine automorphisms defined over a number field $K$ and let $r(f)$ be $D$-ratio of $f$.
        Suppose that $S$ has at least two elements and $r= \displaystyle \max_{l} r(f_l)$. Then, there is a constant $C$ satisfying
        \[
        \sum_{l=1}^k \dfrac{1}{\deg f_l} h\bigl( f(P) \bigr)   > \left( 1+ \dfrac{1}{r}  \right)h(P) - C
        \]
        for all $P \in \af^n (\overline{K})$.
    \end{thm}

    Thus, Silverman's result for preperiodic points \cite[Theorem 4]{S4} is also improved;
     \begin{thm}\label{main2}
        Let $S=\{f_1, \cdots, f_k~|~ f_l :\af^n(K) \rightarrow \af^n(K) \}$ be jointly regular and let $\Phi$ be the monoid of rational maps generated by $S$.
        Define
        \[
        \delta_S : = \left(\dfrac{1}{1+1/r} \right)\sum_{l=1}^k \dfrac{1}{\deg f_l}
        \]
        where $r = \max \bigl(r(f_l)\bigr)$.

        If $\delta_S <1$, then,
        \[
        \Pre(\Phi):= \bigcap_{f \in \Phi} \Pre(f)  ~\text{is a set of bounded height.}
        \]
    \end{thm}

    From now on, we will let $K$ be a number field, let $H$ be an infinity hyperplane of $\af^n$ in the projective space $\pp^n(\overline{K})$ and let $f$ be an affine automorphism unless stated otherwise.

\par\noindent\emph{Acknowledgements}.\enspace
It is a part of my Ph.D. dissertation. The author would like to thank my advisor Joseph H. Silverman for his overall advice.


\section{Preliminaries}

     We need two main ingredients of this paper, the theory of resolution of indeterminacy and the $D$-ratio of rational maps. For details,
     the author refers readers to \cite{Cu} and \cite[II.7]{H} for the resolution of indeterminacy and blowups, and \cite{Le2} for the $D$-ratio.
    \subsection{Blowup and resolution of indeterminacy}
    \begin{thm}[Resolution of Indeterminacy]
        Let $f :V \rightarrow W$ be a rational map between proper varieties such that $V$ is nonsingular. Then there is a proper nonsingular variety $\widetilde{V}$ with a birational morphism $\pi: \widetilde{V} \rightarrow V$ which satisfy that $\phi = f \circ \pi : \widetilde{V} \rightarrow W$ is a morphism.
        \[
                \xymatrix{
                \widetilde{V} \ar[d]_{\pi} \ar[rd]^{\phi} & \\
                V \ar@{-->}[r]_{f} & W
                }
        \]
    \end{thm}

    For notational convenience, we will define the followings;
    \begin{df}
        Let $f :\pp^n \dashrightarrow \pp^n$ be a rational map and let $V$ be a blowup of $\pp^n$ with a birational morphism $\pi :V \rightarrow \pp^n$. We say that a pair $(V,\pi)$ is \emph{a resolution of indeterminacy of $f$} if
        \[ f \circ \pi : V \rightarrow \pp^n\]
        is extended to a morphism $\phi$. And we call the extended morphism $\phi := f \circ \pi$ a \emph{resolved morphism of $f$}.
    \end{df}

    \begin{df}
        Let $pi : W \rightarrow V$ be a birational morphism. We say $\pi$ is a \emph{monoidal trnasformation} if its center scheme is a smooth irreducible subvariety.
    \end{df}

    \begin{thm}[Hironaka]\label{Hironaka}
        Let $f :X \rightarrow Y$ be a rational map between proper varieties such that $V$ is nonsingular. Then, there is a sequence of proper varieties $X_0, \cdots, X_m$ such that
        \begin{enumerate}
        \item $X_0 = X$.
        \item $\rho_i : X_i \rightarrow X_{i-1}$ is a monoidal transformation.
        \item If $T_i$ is the center of blowup of $X_i$, then $\rho_0 \circ \cdots \circ \rho_i (T_i) \subset Z(f)$ on $X$.
        \item $f$ ie extended to a morphism $\widetilde{f}:X_m \rightarrow Y$ on $X_m$.
        \item Consider the composition of all monoidal transformation: $\rho : X_m \rightarrow X$. Then, the underlying set of the center of blowup for $X_m$ is exactly $Z(f)$ on $X$.
        \end{enumerate}
    \end{thm}
    \begin{proof}
        See \cite[Question (E) and Main Theorem II]{Hi}.
    \end{proof}

    \begin{df}
        Let $\pi :V \rightarrow \pp^n$ be a birational morphism. Then, we define $\mathfrak{I}$ is \emph{the center scheme of $\pi$} if
        its corresponding ideal sheaf $\mathcal{S}$ generates $V$:
        \[
        V = \operatorname{Proj}\bigoplus_{d\geq 0}\mathcal{S}^d.
        \]
    \end{df}

    \begin{df}
        Let $\pi : \widetilde{V} \rightarrow V$ be a birational morphism with center scheme $\mathfrak{I}$ and let $D$ be an irreducible divisor of $V$. We define \emph{the proper transformation of $D$ by $\pi$} to be
        \[
        \pi^{\#}D = \overline{\pi^{-1} (D \cap U)}
        \]
        where $U= V - Z \left( \mathfrak{I}\right)$ and $Z \left( \mathfrak{I}\right)$ is the underlying subvariety made by the zero set of the ideal $\mathfrak{I}$.
    \end{df}

    \subsection{$\af^n$-effectiveness and the $D$-ratio}
    \begin{prop}\label{Pic}
        Let $\pi : V \rightarrow \pp^n$ be a birational morphism which is a composition of monoidal transformation. Then, $\Pic(V)$ is a free $\mathbb{Z}$-module. Furthermore, let $H$ be a hyperplane on $\pp^n$ and let $E_i$ be the proper transformation of the exceptional divisor of $i$-th blowup. Then,
        \[
        \{H_V = \pi^\#H, E_1, \cdots, E_r \}
        \]
        is a linearly independent generator of $\Pic(V)$.
    \end{prop}
    \begin{proof}
    \cite[Exer.II.7.9]{H} shows that
    \[
    \Pic(\widetilde{X}) \simeq \Pic(X) \oplus \mathbb{Z}
    \]
    if $\pi:\widetilde{X} \rightarrow X$ is a monoidal transformation. More precisely,
    \[
    \Pic(\widetilde{X}) = \{ \pi^\# D + nE ~|~ D \in \Pic(X)\}.
    \]

    Now suppose that $X_0 = \pp^n$, $\rho_i : X_i \rightarrow X_{i-1}$ is a monoidal transformation. Then, we get the desired result.
    \end{proof}

    \begin{df}
        Let $V$ be a blowup of $\pp^n$, $H$ be a fixed hyperplane of $\pp^n$ and
        \[
        \Pic(V) = \mathbb{Z} H_V \oplus \mathbb{Z}E_1 \oplus \cdots \oplus  \mathbb{Z}E_r.
        \]
        Then, we define the \emph{$\af^n$-effective cone}
        \[
        \operatorname{AFE}(V) = \mathbb{Z}^{\geq 0} H_V \oplus \mathbb{Z}^{\geq 0} E_1 \oplus \cdots \oplus  \mathbb{Z}^{\geq 0} E_r
        \]
        where $\mathbb{Z}^{\geq 0}$ is the set of nonnegative integers.
        We say a divisor $D$ of $V$ is \emph{$\af^n$-effective}
        if $D\in \operatorname{AFE}(V)$ and denote it by
        \[
        D \succ 0.
        \]
        Moreover, we will say
        \[
        D_1 \succ D_2
        \]
        if $D_1 - D_2$ is $\af^n$-effective.
    \end{df}

        \begin{prop}\label{af effec prop} Let $V$ be a blowup of $\pp^n$ with birational morphism $\pi:V \rightarrow \pp^n$ and $D, D_i \in \Pic(V)$.\\
        {\rm(1)  (Effectiveness)} If $D$ is $\af^n$-effective, then $D$ is effective.\\
        {\rm(2)  (Boundedness)} If $D$ is $\af^n$-effective, then
        $h_D(P)$ is bounded below on $V \setminus \left( H_V \cup \left( \bigcup_{i=1}^r E_i \right) \right) $. \\
        {\rm(3) (Transitivity)} If $D_1 \succ D_2$ and $D_2 \succ D_3$ , then $D_1 \succ D_3$\\
        {\rm (4) (Funtoriality)} If $W$ is a blowup of $V$, a map $\rho:W \rightarrow V$ is a birational morphism and $D_1 \succ D_2$, then
        $\rho^*D_1 \succ \rho^* D_2.$
    \end{prop}
    \begin{proof}
        See \cite[Proposition 3.3]{Le2}
    \end{proof}

   \begin{df}
        Let $f :\pp^n \dashrightarrow \pp^n$ be a rational map with $Z(f) \subset H$, let $(V, \pi_V)$ be a resolution of indeterminacy of $f$ and let $\phi_V$ is a resolved morphism.
            \[
                \xymatrix{
                V \ar[d]_{\pi_V} \ar[rd]^{\phi_V}\\
                \pp^n \ar@{-->}[r]_{f}  & \pp^n }
            \]
        Suppose that
            \[
            \pi_V^*H = a_0 H_V + \sum_{i=1}^r a_i E_i \quad \text{and} \quad \phi_V^*H = b_0 H_V + \sum_{i=1}^r b_i E_i
            \]
        where $a_i,b_i$ are nonnegative integers. If all $b_i$ are nonzero for all $i$ satisfying $a_i \neq 0$, we define \emph{the $D$-ratio} of $\phi_V$,
        \[r(\phi_V)= {\deg \phi_V} \cdot \max_i \left( \dfrac{a_i}{ b_i} \right).\]
        Otherwise; if there is an $i$ satisfying $a_0\neq 0$ and $b_i=0$, define
        \[r(\phi_V) = \infty.\]
    \end{df}

    \begin{lem}\label{invariant}
        Let $(V,\pi_V)$ and $(W,\pi_V)$ be resolutions of indeterminacy with resolved morphisms $\phi_V = f \circ \pi_V$ and $\phi_W = f \circ \pi_W$ respectively.
            \[
                \xymatrix{
                W \ar[d]_{\pi_W}  \ar[rd]^{\phi_W} & &  V \ar[d]^{\pi_V} \ar[ld]_{\phi_V} \\
                \pp^n \ar@{-->}[r]_{f} & \pp^n   & \pp^n \ar@{-->}[l]^{f}
                }
            \]
        Then,
    \[
    r(\phi_V)= r(\phi_W).
    \]
    \end{lem}
    \begin{proof}
    See \cite[Lemma 4.3]{Le2}
    \end{proof}

    \begin{df}\label{Dratio}
        Let $f :\pp^n \dashrightarrow \pp^n$ be a rational map with $Z(f) \subset H$. Then, we define \emph{the $D$-ratio of $f$}, \[r(f)= r(\phi_V)\]
        for any resolution of indeterminacy $(V,\pi_V)$ of $f$ with resolved morphism $\phi_V$.
    \end{df}

    \begin{prop}
       Let $f,g:\pp^n \dashrightarrow \pp^n$ be rational maps with $Z(f),Z(g) \subset H$.  Then,
        \begin{enumerate}
            \item $r(f)=1$ if $f$ is a morphism.
            \item $r(f) \in [1,\infty] $ .
            \item  $\dfrac{r(f)}{\deg f} \cdot \dfrac{r(g)}{\deg g} \geq \dfrac{r(g\circ f)}{\deg (g\circ f) }$.
            \item If $g$ is a morphism and $f$ is a rational map on $\pp^n$, then
                $   r(g\circ f)= r(f)$.
         \end{enumerate}
    \end{prop}

    \begin{ex}
        Let $f :\af^n \rightarrow \af^n$ be an affine automorphism with the inverse map $f^{-1}:\af^n \rightarrow \af^n$. Then,
        $r(f) = \deg f \times \deg f^{-1}.$ (For details, see \cite{Le}.) For example, a H\'{e}non map
        \[
        f_H(x,y,z) = (z, x+z^2 ,y +x^2)
        \]
        is an example of regular affine automorphism with the inverse map
        \[
        f_H^{-1}(x,y,z) = (y-x^2, z-(y-x^2)^2 ,x).
        \]
        Thus,
        \[
        r(f_H) = r(f_H^{-1}) = \deg f_H \times \deg f_H^{-1} = 2 \times 4 = 8.
        \]
    \end{ex}

                \begin{ex}
                    Let $f[x,y,z]= [x^2,yz,z^2]$. Then, the indeterminacy locus is $P=[0,1,0]$. Then, the blowup $V$ along closed scheme corresponding ideal sheaf $(yz,x^2)$ will resolves
                    indeterminacy, which is a successive blowup along $P$ and $H^\# \cap E_1$.
                    $$f_1[x,y,z][x_1, z_1]= [x_1x,z_1y,z_1z]$$
                    $$\phi = f_2[x,y,z][x_1, z_1][x_2,z_2]= [x_2,z_2y,x_2z_1^2]$$

                    Let $E_1, E_2$ be the exceptional divisors on each step.
                    
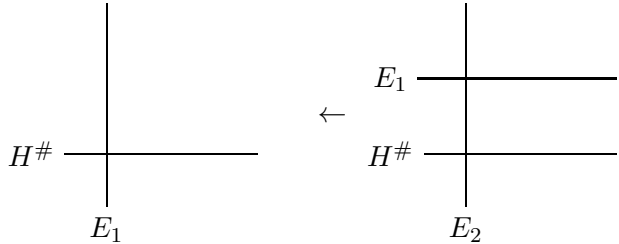
\[
                        \xy <1cm,0cm>:
                            (0,1) *+{H^\#}; (3,1) **@{-},
                            (1,0) *+{E_1}; (1,3) **@{-},
                            (4,1.5)*+{\leftarrow}
                        \endxy
                        \xy <1cm,0cm>:
                            (0,1) *+{H^\#}; (3,1) **@{-},
                            (1,0) *+{E_2}; (1,3) **@{-},
                            (0,2) *+{E_1}; (3,2) **@{-}
                        \endxy
                    \]

                    Then, the intersection number $E_2^2 = -1$, $E_1^2=-2$, $(H^\#)^2=-1$, $H^\# \cdot E_1=0$ and
                    $H^\#\cdot E_2 = E_1 \cdot E_2 = 1$

                    Furthermore,
                    $$H^\# \cdot \phi^*H = \phi_*H^\# \cdot H = 0,$$
                    $$E_1 \cdot \phi^*H = \phi_*E_1 \cdot H = 0,$$
                    $$E_2 \cdot \phi^*H = \phi_*E_2 \cdot H = 1 . $$

                    Since $\Pic (V) = \langle H^\# , E_1, E_2 \rangle$, we may assume that
                    $$\phi^*H = aH^\# + bE_1 + cE_2$$
                    Then, by previous facts,
                    $$\phi^*H\cdot H^\# = -a+c = 0, \quad \phi^*H\cdot E_1 = a - 2b = 0 .$$
                    Therefore,
                    $$\phi^*H = 2H^\# + E_1 + 2E_2, \quad \pi^*H = H^\# + E_1 + 2E_2$$
                    and hence
                    $$r(f) = 2 \times 1 = 2$$

                \end{ex}


\section{Jointly Regular Families of Rational maps}

     \begin{proof}[Proof of Theorem~\ref{main}]

     For notational convenience, let
     \begin{itemize}
     \item $d_l = \deg f_l$
     \item $r_l = r(f_l)$
     \item $(V_l, \pi_l)$ be a resolution of indeterminacy of $f_l$ constructed by Theorem~\ref{Hironaka}; assume $\pi_l$ is a composition of monoidal transformation and $\{\pi_l^\# H = H_{V_l}, E_{l1}, \cdots, E_{ls_l}\}$ is the generator of $\Pic(V_l)$ given by Proposition~\ref{Pic}.
     \item $\phi_l$ be the resolved morphism of $f_l$ on $V_l$.
     \item
        \[
        \pi_l^*H = a_0H_{V_l} + \sum_{i=1}^{s_l} a_{li} E_{li} \quad \text{and} \quad \phi_l^*H = b_0 H_{V_l} + \sum_{i=1}^{s_l} b_{li} E_{li}
        \]
        in $\Pic(V_l)= \mathbb{Z} \pi_{l}^\#H \oplus \mathbb{Z} E_{l1} \oplus \cdots \oplus \mathbb{Z} E_{ls_l}$.
     \end{itemize}

         We can easily check that $a_0 =1$ and $b_0=d_l$ from ${\pi_l}_*\pi_l^*H = H$ and ${\pi_l}_*\phi_l^*H = \deg \phi_l \cdot H$ For details, see \cite[Proposition 4.5.(2)]{Le2}.

        Let $T_l$ be the center scheme of blowup for $V_l$ and $W$ is the blowup of $\pp^n$ whose center scheme is $\sum T_l$. Then, $W$ is a blowup of $V_l$ for all $l$. Furthermore, since the underlying set of $T_l$ is exactly $Z(f_l)$, the underlying set of $\sum T_l = \cup Z(f_l)$. Let $\rho_l : W \rightarrow V_l$, $\pi_W$ be the monoidal transformations defined by construction of $W$:
        \[
            \xymatrix{
                & &  W \ar[dd]|{\pi_{W}}
                \ar[ld]_{\rho_l} \ar[rd]^{\rho_{l'}} \ar@/_2pc/[lldd]_{\widetilde{\phi_l}} \ar@/^2pc/[rrdd]^{\widetilde{\phi_{l'}}}& & \\
                & V_l \ar[ld]_{\phi_l} \ar[rd]^{\pi_l} & & V_{l'} \ar[rd]^{\phi_{l'}} \ar[ld]_{\pi_{l'}} & \\
                \pp^n & & \pp^n \ar@{-->}[ll]^{f_l} \ar@{-->}[rr]_{f_{l'}} & & \pp^n
            }
        \]

        Then, still $W$ is a blowup of $\pp^n$ and hence $\Pic(W)$ is generated by $\pi_W$ and the irreducible compoenets of the exceptional divisor:
        \[
        \Pic(W) = \mathbb{Z} \pi_{W}^\#H \oplus \mathbb{Z} F_1 \oplus \cdots \oplus \mathbb{Z} F_s
         \]
         where $F_j$ are irreducible components of exceptional divisor of $W$. Thus, we can represent $\pi_W^* H$ as follows:
         \[
         \pi_W^*H = \pi_{W}^\#H + \sum_{j=1}^{s} \alpha_{j} F_{j}.
        \]

        To describe $\phi_l^*H$ precisely, let's define sets of indices
        \[
        \mathcal{I}_l = \{1 \leq j \leq s ~|~ \pi_W(F_j) \subset Z(f_l) \} \quad \text{and} \quad
        \mathcal{I}^c_l = \{1 \leq j \leq s ~|~ \pi_W(F_j) \not \subset Z(f_l) \}.
        \]
        By definition, it is clear that
        \[
        \mathcal{I}_{l} \cup \mathcal{I}_{l}^c = \{ 1,\cdots, s\} \quad \text{and}\quad \mathcal{I}_{l} \cap \mathcal{I}_{l}^c = \emptyset.
        \]
        Thus, we can say
        \[
        \widetilde{\phi}_l^*H = d_l \pi_{W}^\#H + \sum_{j=1}^{s} \beta_{lj} F_j
        = d_l \pi_{W}^\#H + \sum_{j \in \mathcal{I}^c_l} \beta_{lj} F_{j} +\sum_{j \in \mathcal{I}_l} \beta_{lj} F_{j}.
        \]

        Moreover, we have the following lemmas;
        \begin{lem}\label{indices}
        \[
        \bigcup_{l=1}^k \mathcal{I}_{l} = \bigcup_{l=1}^k \mathcal{I}^c_l = \{ 1, \cdots , s\}.
        \]
        \end{lem}
        \begin{proof}
        $\bigcup_l \mathcal{I}_l =  \{ 1, \cdots , s\}$ is clear; because the underlying set of the center scheme of $W$ is $\cup Z(f_l)$,         $\cup \pi_W(F_j) = \pi_W(\cup F_j) = \cup Z(f_l)$.

        Suppose $\bigcup_l \mathcal{I}^c_l \subsetneq \{ 1, \cdots , s\}$. Then, there is an index $k_0$ satisfying
        $\pi_W(F_{k_0}) \subset Z(f_{l})$ for all $l$. This implies $ \pi_W(F_{k_0}) \subset Z(f_{l})$ for all $l$ and hence
        $\emptyset \neq \pi_W(F_{k_0}) \subset \bigcap_{l} Z(f_l) $
        which contradicts to that $S$ is jointly regular.
        \end{proof}

        \begin{lem}\label{coefficients}
        Let $\alpha_j$ and $\beta_{lj}$ be the coefficients of $F_j$ in $\pi_V^*H$ and $\widetilde{\phi}_l^*H$ respectively. Then,
        \[
        d_l \dfrac{\alpha_j}{\beta_{lj}} \leq r_l.
        \]
        Especially, if $j \in \mathcal{I}^c_l$, then
        \[d_l \alpha_j = \beta_{lj}.\]
        \end{lem}
        \begin{proof}
        By definition of the $D$-ratio, the first inequality is clear:
        \[
        r_l = d_l \cdot \max_i \left( \dfrac{\alpha_i}{\beta_{li}} \right) \geq d_l \cdot  \dfrac{\alpha_j}{\beta_{lj}}.
        \]

        Now, suppose that
        \[
        \begin{array}{llclc}
        \rho_l^* \pi_l^\#H &=& \gamma_{l00}\pi_W^\# H + \sum_{j=1}^s \gamma_{l0j} F_j &=& \gamma_{l00}\pi_W^\# H + \sum_{j \in \mathcal{I}^c_l} \gamma_{l0j} F_j
        +\sum_{j \in \mathcal{I}_l} \gamma_{l0j}F_j \\  
        &&&&\\      
        \rho_l^* E_{li} &=& \gamma_{li0}\pi_W^\#H + \sum_{j=1}^s  \gamma_{lij}F_j &=&  \gamma_{li0}\pi_W^\#H + \sum_{j \in \mathcal{I}^c_l} \gamma_{lij}F_j+  \sum_{j \in \mathcal{I}_l} \gamma_{lij}F_j.
        \end{array}
        \]
        First of all, $\gamma_{l00}=1$ and $\gamma_{li0}=0$ for all $i \neq 0$; 
        if $i \neq 0$, ${\pi_W}_*\left( \rho_l^*E_i \right) =0$ because $\pi_W\bigl( \rho_l^*E_i \bigr) \subset \cup Z(f_l)$. On the other hand,
        \[
        \pi_W^\# \left(\gamma_{li0}\pi_W^\#H + \sum_{j=1}^s  \gamma_{lij}F_j \right) = \gamma_{li0} H.
        \]
        Hence, $\gamma_{li0}=0$. For $\gamma_{l00}$, we have
        \[
        {\pi_W}_* \bigl( \pi_W^*H \bigr) = H
        \]
        because $\pi_W$ is one-to-one outside of the center of blowup of $W$. Therefore, 
        \[
        {\pi_W}_* \bigl( \rho_l^* \pi_l^\#H \bigr) = {\pi_W}_* \left( \pi_W^*H - \sum_{j=1}^s a_{li} \rho_l^*E_{li} \right) = H
        \]
        and hence $\gamma_{l00}=1$.
        
        Moreover, because $\pi_l(E_{li}) \subset Z(f_l)$ and $\pi_W(F_j) \not \subset Z(f_l)$ for any $j \in \mathcal{I}^c_l$, the multiplicity of $\rho_l(F_j)$ on $E_l$ is zero and hence $\gamma_{lij} =0$. Thus, we can say
        \[
        \rho_l^* E_{li} =  \sum_{j \in \mathcal{I}_l} \gamma_{lij}F_j.
        \]

        Since $\widetilde{\phi}_l = \rho_l \circ \phi_l$ and $\pi_W = \rho_l \circ \pi_l$, we have
        \[
        \widetilde{\phi}_l^*H = \rho_l^* \phi_l^*H = \rho_l^*\left(d_l  \pi_l^\#H + \sum_{i=1}^{s_l} b_{li} E_{li} \right)
        = d_l \rho_l^*\pi_l^\#H + \sum_{i=1}^{s_l} b_{li} \rho_l^* E_{li}.
        \]
        Thus,
        \begin{eqnarray*}
        \displaystyle \pi_W^* H & = & \rho_l^*\pi_l^* H  \\
        & = & \rho_l^* \left( \pi_l^\# H + \sum_{i=1}^{s_l} a_{li} E_{li} \right) \\
        & = & \left( \pi_W^\# H + \sum_{j \in \mathcal{I}_l} \gamma_{l0j} F_j + \sum_{j \in \mathcal{I}^c_l} \gamma_{l0j} F_j \right)
        + \sum_{i=1}^{s_l}  a_{li}\left( \sum_{j \in \mathcal{I}_l} \gamma_{lij}F_j \right) \\
        & = & \pi_W^\# H + \sum_{j \in \mathcal{I}^c_l} \gamma_{l0j} F_j + \sum_{j \in \mathcal{I}_l}\left(  \sum_{i=0}^{s_l}  a_{li} \gamma_{lij} \right) F_j,
        \end{eqnarray*}
        and
        \begin{eqnarray*}        
        \widetilde{\phi}_l^*H &=& \rho_l^*\phi_l^* H \\
        & = & \rho_l^* \left( d_l \pi_l^\# H + \sum_{i=1}^{s_l} b_{li} E_{li} \right) \\
        & = & d_l \left( \pi_W^\# H + \sum_{j \in \mathcal{I}_l} \gamma_{l0j} F_j + \sum_{j \in \mathcal{I}^c_l} \gamma_{l0j} F_j \right)
        + \sum_{i=1}^{s_l}  b_{li} \left( \sum_{j \in \mathcal{I}_l} \gamma_{lij}F_j \right) \\
        & = & d_l \pi_W^\# H + \sum_{j \in \mathcal{I}^c_l} d_l\gamma_{l0j} F_j + \sum_{j \in \mathcal{I}_l}\left(\sum_{i=0}^{s_l}  b_{li}
        \gamma_{lij} \right) F_j.
        \end{eqnarray*}
        
        Therefore, 
        \[
        d_l \alpha_j  = d_l \sum_{j \in \mathcal{I}^c_l} \gamma_{l0j}  = \beta_j \quad \text{for all}~j \in \mathcal{I}^c_l.
        \]
        \end{proof}
        We now complete the proof of Theorem~\ref{main}. Let $ r = \max r_l$. Note that 
        \[
        p_0 \pi_W^\#H +\displaystyle \sum_{j=1}^s p_j F_j \succ q_0 \pi_W^\#H+ \displaystyle \sum_{j=1}^s q_j F_j
         \]
         if $p_j \geq q_j$ for all $j=0, \cdots, s$. Then, we have
        \[
        \begin{array}{lll}
        \sum_{l=1}^k \dfrac{1}{d_l}\widetilde{\phi}_l^*H
                          & =   \displaystyle \sum_{l=1}^k \left[ \pi_W^\#H
                          +  \sum_{j \in \mathcal{I}^c_l}  \left(\dfrac{\beta_{lj}}{d_l} F_{j}\right)
                          +   \sum_{j \in \mathcal{I}_l} \left(\dfrac{\beta_{lj}}{d_l} F_{j}\right) \right] &  \\
                          & \succ  \displaystyle \sum_{l=1}^k \pi_W^\#H
                          + \sum_{l=1}^k  \sum_{k \in \mathcal{I}^c_l}  \alpha_{j} F_{j}
                          + \sum_{l=1}^k \left( \sum_{j \in \mathcal{I}_l} \dfrac{ \alpha_{j} }{r_l} F_{j}\right) & (\because ~ Lemma~\ref{coefficients})\\
                          & \succ \displaystyle k \pi_W^\#H
                          + \sum_{l=1}^k \sum_{j \in \mathcal{I}^c_l}  \alpha_{j} F_{j}
                          +   \sum_{l=1}^k  \left(  \sum_{j \in \mathcal{I}_l}   \dfrac{ \alpha_{j} }{r} F_{j}\right) &  (\because ~ r \geq r_l)\\
                          & \succ \displaystyle k \pi_W^\#H
                          + \sum_{j=1}^s  \alpha_{j} F_{j}
                          +   \dfrac{ 1 }{r}   \sum_{j=1}^s   \alpha_{j}  F_{j} & (\because ~ Lemma~\ref{indices}) \\
                          & \succ \displaystyle
                          \left( 1+  \sum \dfrac{ 1 }{r} \right)  \pi_W^* H &
        \end{array}
        \]
        and hence
        \[
        D = \sum_{l=1}^k \dfrac{1}{d_l}\widetilde{\phi}_l^*H -\left( 1+  \sum \min \dfrac{ 1 }{r_l} \right)  \pi_W^* H
        \]
        is an $\af^n$-effective divisor.

        Thus, by Proposition~\ref{af effec prop}, $h_D$ is bounded below on $\pi_W^{-1} \af^n$. Therefore, there is a constant $C$ such that
        \begin{eqnarray*}
        h_D(Q) &=& \sum_{l=1}^k \dfrac{1}{d_l} h_{\widetilde{\phi}_l^*H}(Q) - \left( 1+  \sum \min \dfrac{ 1 }{r_l} \right)  h_{\pi_W^* H}(Q) \\
        &=& \sum_{l=1}^k \dfrac{1}{d_l} h_{^*H}\bigl( \widetilde{\phi}_l(Q) \bigr)
        - \left( 1+  \sum \min \dfrac{ 1 }{r_l} \right)  h_{ H}\bigl( \pi_W^*(Q)\bigr) \\
        &>& C
        \end{eqnarray*}
        for all $Q \in \pi_W^{-1}(\af^n)(\overline{K})$.  Finally, for $P = \pi(Q)$, $\widetilde{\phi}_l(Q) = f(P)$ and $\pi_W(Q) = P$ and hence
        \[
        \sum_{l=1}^k \dfrac{1}{d_l} h_H(P) - \left( 1+  \sum \min \dfrac{ 1 }{r_l} \right)  h_H(P) > C.
        \]
        \end{proof}

    \begin{ex}
        Let 
        \[
        f_1= (z,y+z^2, x+(y+z^2)^2), \quad f_2 = (x,y^2,z), \quad f_3 = (x^2+y,x+y,z^3).
        \]
        Then, the $r(f_1) = 8$, $r(f_2) =2$ and $r(f_3) = 3/2$. (For details of the $D$-ratio calculation, see \cite{Le2}.) Therefore, 
        \[
        h\bigl( (z,y+z^2, x+(y+z^2)^2) \bigr) + h\bigl( (x,y^2,z) \bigr) + h\bigl( (x^2+y,x+y,z^3) \bigr) \geq 
        \left(1 + \dfrac{1}{8} \right) h\bigl( (x,y,z) \bigr) -C
        \]
        for some constant $C$.
    \end{ex}

    \begin{cor}\label{liminf}
        Let $S$ be a jointly regular set of affine morphisms. Then,
        \[
        \kappa (S) :=
        \liminf_{\substack{P\in \af^n(\overline{\mathbb{Q}})\\h(P)\rightarrow \infty}}
        \sum_{f \in S}\dfrac{1}{\deg f} \dfrac{ h\bigl( f(P) \bigr)}{h(P)} \geq 1+ \dfrac{1}{r}
        \]
        where $r = \displaystyle \max_{f\in S} r(f)$.
    \end{cor}

    \begin{rem}
        Corollary~\ref{liminf} may not be the exact limit infimum valus. For example, If there is a subset $S' \subset S$
        such that $S'$ is still jointly regular and $\displaystyle \max_{f\in S'} r(f) < \displaystyle \max_{f\in S} r(f)$, then
        \[
        \kappa(S) \geq \kappa(S')        \geq 1 + \dfrac{1}{r'} > 1+ \dfrac{1}{r}.
        \]
     \end{rem}
     
     \begin{ex}
        We have some examples for $\kappa(S) = 1 + \displaystyle \min_{f\in S}\left(\dfrac{1}{r(f)} \right)$;
        \begin{enumerate}
        \item $S=\{f,g\}$ where $f,g$ are morphisms. 
        
                If $f,g$ are morphism, then $r(f)=r(g)=1$. Therefore, 
                \[
                \dfrac{1}{\deg f}h\bigl(f(P) \bigr) + \dfrac{1}{\deg g}h\bigl(g(P) \bigr) = h(P) + h(P) + O(1).
                \]
                
        \item $S= \{f,f^{-1} \}$ where $f$ is a regular affine automorphism and $f^{-1}$ is the inverse of $f$.
        
                It is proved by Kawaguchi. See \cite{K2}.
                
        \end{enumerate}
    \end{ex}

\section{An application to arithmetic dynamics}

    This result is a generalization of \cite[Section 4]{S4}. The only difference is that we have an improved inequality for jointly regular family. The proof is almost the same.

    Fix an integer $m\geq 1$ and let $S=\{f_1, \cdots, f_k\} \subset \Rat^n(H)$ be a jointly regular family. For each $m\geq 0$, let $W_m$ be the collection of ordered $m$-tuples chosen from $\{1, \cdots, k\}$,
    \[
    W_m = \bigl\{ (i_1, \cdots, i_m) ~|~ i_j \in \{1, \cdots, k\} \bigr\} = \bigl\{1, \cdots, k \bigr\}^m,
    \]
    and let
    \[
    W_* = \bigcup_{m\geq 0} W_m.
    \]
    Thus $W_*$ is the collection of words on $r$ symbols.

    For any $I=(i_1, \cdots, i_m) \in W_m$, let $f_I$ denote the corresponding composition of the rational amps $f_1, \cdots f_k$,
    \[
    f_I = f_{i_1} \circ \cdots \circ f_{i_m}.
    \]
    \begin{df}
        We denote \emph{the monoid of rational maps generated by $f_1, \cdots, f_k$ under composition} by
        \[
        \Phi = \{\phi = f_I ~|~ I \in W_*\}.
        \]
        Let $P\in \af^n$. \emph{The $\Phi$-orbit of $P$} is
        \[
        \Phi(P) = \{ \phi(P) ~|~ \phi \in \Phi \}.
        \]
        \emph{The set of (strongly) $\Phi$-preperiodic points }is the set
        \[
        \Pre(\Phi) = \{P \in \af^n ~|~ \Phi(P)~\text{is finite}\}.
        \]
    \end{df}

    \begin{proof}[Proof of Theorem~\ref{main2}]
     By Theorem~\ref{main}, we have a constant $C$ such that
    \begin{eqnarray}\label{ineq1}
    0 \leq  \left(\dfrac{1}{1+\frac{1}{r}}\right) \sum_{l=1}^k \dfrac{1}{d_l} h\bigl( f_l(Q) \bigr) - h(Q) +C \quad \text{for all}~Q\in \af^n.
    \end{eqnarray}
    Note that if $r=\infty$, then $\left(\dfrac{1}{1+\frac{1}{r}}\right) = 1$ and theorem holds because of \cite[Section 4]{S4}. Thus, we may assume that $r$ is finite.

    We define a map $\mu:W_* \rightarrow \mathbb{Q}$ by the following rule:
    \[
    \mu_I = \mu_{(i_1, \cdots, i_m)} = \prod d_l^{p_{I,l}}
    \]
    where $p_{I,l} = - |\{t~|~ i_t = l\}|$. Then, by definition of $\delta$ and $\mu_I$, the following is true:
    \[
    \delta^m = \left[ \left(\dfrac{r}{r+1} \right)\sum_{l=1}^k \dfrac{1}{d_l} \right]^m
    = \left(\dfrac{r}{r+1} \right)^m \sum_{I\in W_m} \dfrac{1}{\deg f_{i_1} \cdots \deg f_{i_m}}
    = \left(\dfrac{r}{r+1} \right)^m \sum_{I\in W_m} \mu_I.
    \]

    Let $P\in \af^n(\overline{\mathbb{Q}})$. Then, (\ref{ineq1}) holds for $f_I(P)$ for all $I\in W_m$:
    \begin{eqnarray*}
    0 \leq  \left(\dfrac{r}{r+1}\right) \sum_{l=1}^k \dfrac{1}{d_l} h\bigl( f_l(f_I(P) ) \bigr) - h(f_I(P)) +C.
    \end{eqnarray*}
    and hence
    \begin{eqnarray}\label{mainineq}
    0 \leq  \sum_{m=0}^M \sum_{I\in W_m}\mu_I \left(\dfrac{r}{r+1}\right)^m \left[ \sum_{l=1}^k \dfrac{1}{d_l} h\bigl( f_l(f_I(P) ) \bigr) -
    \left(1+ \dfrac{1}{r}\right) h(f_I(P)) +C \right].
    \end{eqnarray}

    The main difficulty of the inequality is to figure out the constant term. From the definition of $\delta$, we have
    \begin{eqnarray*}
    \sum_{m=0}^{M-1} \left(\dfrac{r}{r+1}\right)^m \sum_{I\in W_m} \mu_I = \sum_{m=1}^M \delta^m \leq \dfrac{1}{1-\delta}.
    \end{eqnarray*}
    Now, do the telescoping sum and most terms in (\ref{mainineq}) will be canceled;
    \begin{eqnarray*}\label{identity}
     &&\left(\sum_{m=0}^{M-1} \sum_{I \in W_m} \left(\dfrac{r}{r+1}\right)^m\mu_I \sum_{l=1}^k \dfrac{1}{d_k} h\bigl( f_lf_I(P) \bigr)\right)
     - \left(\sum_{m=1}^{M} \sum_{I \in W_m} \left(\dfrac{r}{r+1}\right)^{m-1}\mu_I h\bigl( f_I(P) \bigr) \right)\\
     &&= \left(\sum_{m=0}^{M-1} \sum_{I \in W_m} \left(\dfrac{r}{r+1}\right)^m \sum_{l=1}^k \dfrac{\mu_I}{d_l} h\bigl( f_lf_I(P) \bigr)\right)
     - \left(\sum_{m=0}^{M-1} \sum_{I \in W_m} \sum_{l=1}^k \left(\dfrac{r}{r+1}\right)^m \dfrac{\mu_{I}}{d_l} h\bigl( f_lf_I(P) \bigr)\right) \\
     &&=0
    \end{eqnarray*}
    Therefore, the remaining terms in (\ref{mainineq}) are
    \begin{eqnarray*}
    0 &\leq&  \left[ \sum_{I\in W_M} \left(\dfrac{r}{r+1}\right)^M \mu_I \sum_{l=1}^K \dfrac{1}{d_l} h\bigl( f_l(f_I(P) ) \bigr)
    \right] - h(P) +  \sum_{I\in W_M} \left(\dfrac{r}{r+1}\right)^M \mu_I C \\
    &\leq&  \left[ \sum_{I\in W_M} \left(\dfrac{r}{r+1}\right)^M \mu_I \sum_{l=1}^k \dfrac{1}{d_l} h\bigl( f_l(f_I(P) ) \bigr)
    \right] - h(P) + \dfrac{1}{1-\delta} C,
    \end{eqnarray*}
   \[
   \sum_{I\in W_M} \left(\dfrac{r}{r+1}\right)^M \mu_I \sum_{l=1}^k \dfrac{1}{d_l} = \left(\dfrac{r}{r+1}\right)^M \sum_{I\in W_{M+1}} \mu_I
   = \left(1+\dfrac{1}{r}\right) \delta^{M+1}.
   \]
    Define the height of the images of $P$ by the monoid $\Phi$:
   \[
   h(\Phi(P)) = \sup_{R\in \Phi(P)} h(R).
   \]
    Then, if $P\in \Pre(\Phi)$, $h(\Phi(P))$ is finite and hence we have an upper bound for $h(P)$:
    \begin{eqnarray*}
    h(P) &\leq&  \left[ \sum_{I\in W_M} \left(\dfrac{r}{r+1}\right)^M \mu_I \sum_{l=1}^k \dfrac{1}{d_l}\right]  h\bigl( \Phi(P) \bigr)+ \dfrac{1}{1-\delta} C\\
    &\leq&  \left(1+\dfrac{1}{r}\right) \delta^{M+1}  h\bigl( \Phi(P) \bigr)+ \dfrac{1}{1-\delta} C.
    \end{eqnarray*}
    By assumption, $\delta<1$ and $h \bigl(\Phi(P)\bigr)$ is finite, so letting $M \rightarrow \infty$ shows that
    $h(P)$ is bounded by a constant that depends only on $S$.
    \end{proof}

\end{document}